\documentclass[12pt]{amsart}

\newcommand{\NN}{\mathbb{N}}
\newcommand{\VV}{\mathcal{V}}

\newcommand{\pp}{\mathbb{P}}
\newcommand{\RR}{\mathbb{R}}
\newcommand{\GG}{\mathcal{G}}
\newcommand{\ZZ}{\mathbb{Z}}

\newcommand{\sop}[1]{\mbox{{\em supp\/}}(#1)}

\newcommand{\es}[1]{\mbox{{\em ess\,lim\/}}#1}
\newcommand{\elimsup}[1]{\mbox{{\rm ess\,lim\,sup\/}}#1}
\newcommand{\eliminf}[1]{\mbox{{\rm ess\,lim\,inf\/}}#1}
\newcommand{\esup}[1]{\mbox{{\rm ess\,sup\/}}#1}

\newtheorem{definition}{\sc Definition}[section]
\newtheorem{teo}{\sc Theorem}[section]
\newtheorem{prop}{\sc Proposition}[section]

\begin{document}

\title{On Hilbert extensions of Weierstrass' theorem with weights}

\author{YAMILET QUINTANA}
\thanks{{\rm Research partially supported by  DID-USB under  Grant DI-CB-015-04}}
\email{yquintana@usb.ve}
\begin{abstract}
In this paper we study the set of functions $\GG$-valued which can
be approximated by $\GG$-valued continuous functions in the norm
$L^\infty_{\GG}(I,w)$, where $I$ is a compact interval, $\GG$ is a
real and separable Hilbert space and $w$ is certain  $\GG$-valued
weakly measurable weight. Thus, we obtain a new extension of
celebrated Weierstrass approximation theorem.
\\

\hspace{-.6cm} {\it Key words and phrases.} Weierstrass' theorem,
$\GG$-valued weights, $\GG$-valued polynomials, $\GG$-valued continuous functions.\\

\hspace{-.6cm} {\it 2001 Mathematics Subject Classification.}
Primary 41, 41A10. Secondary 43A32, 47A56.
\end{abstract}

\maketitle \markboth{YAMILET QUINTANA  \mbox{} }{ON HILBERT
EXTENSIONS OF WEIERSTRASS' THEOREM WITH WEIGHTS}

\rightline{{\it To Ma\'{\i}ta, in memoriam.}}

\section{Introduction.}

If $I$ is any compact interval, Weierstrass' approximation theorem
says that $C(I)$ is the largest set of functions which can be
approximated by polynomials in the norm $L^\infty(I)$, if we
identify, as usual, functions which are equal almost everywhere.
Weierstrass proved this theorem in 1885. Also, in that time he
proved the density of trigonometric polynomials in the class of
$2\pi$-periodic continuous real-valued functions. These results were
-in a sense- a counterbalance to Weierstrass' famous example of 1861
on the existence of a continuous nowhere differentiable function
(see \cite{Bai}).

Two subject of interest for Weierstrass were complex function theory
and the possibility of representing functions by power series. The
result obtained in his paper in 1885 should be viewed from that
perspective, moreover, the title of the paper  emphasizes such
viewpoint (his paper was titled {\it On the possibility of giving an
analytic representation to an arbitrary function of real variable},
see \cite{L}). Weierstrass' perception on analytic functions was
`functions that could be represented by power series'. This paper of
Weierstrass was reprinted in Weierstrass' Mathematische Werke
(collected works) with some notable additions, for example a short
`introduction'. This reprint appeared in 1903 and contains the
following statement:

{\it The main result of this paper, restrict to the one variable
case, can be summarized as follows:}

{ \it Let $f\in C(\RR)$. Then there exists a sequence $f_1 ,f_2
,\ldots $ of entire functions for which}
$$f(x)=\sum_{i=1}^{\infty}f_{i}(x),$$
{\it for each $x\in\RR$. In addition the convergence of above sum is
uniform on every finite interval.}

Let us observe that there is not mention of fact that the $f_i$ may
be assumed to be polynomial. We state Weierstrass' approximation
theorem, not as given in his paper, but as it is currently stated
and understood.

\begin{teo}
\label{t1}
(K. Weierstrass).\\
Given  $f:[a,b]\rightarrow \RR$ continuous and an arbitrary
$\epsilon>0$ there exists an algebraic polynomial $p$ such that
\begin{equation}
\left|f(x)-p(x)\right|\leq \epsilon, \,\, \forall\,\, x\in [a,b].
\end{equation}
\end{teo}

Two papers of Runge published about the same time also provided a
proof of this result. But unfortunately, the theorem was not titled
Weierstrass-Runge theorem. The impact of Weierstrass' approximation
theorem on the mathematical world was immediate: there were later
proofs of famous mathematicians such as Picard (1891), Volterra
(1897), Lebesgue (1898), Mittag-Leffler (1900), Landau (1908), de la
Valle\'e Poussin (1912). The proofs more commonly taken at level of
undergraduate courses in Mathematics are those of Fej\'er (1900) and
Bernstein (1912) (see, for example \cite{Bar}, \cite{Ch}, \cite{L} ,
\cite{P} or \cite{Q}).

There have been many improvements, generalizations and ramifications
of Weierstrass' approximation theorem. For instance, if $f$ is a
$\VV$-valued function, with $\VV$ a real (or complex) finite
dimensional linear space, and if $f$ is a function of several real
variables.  Each one of these cases have an easy formulation of
results. While,  the case in which  $f$ is a function of several
complex variables requires a more profound study, with skillful
adaptations of both hypothesis and conclusion. A detailed
presentation of such results  may be found in \cite{S1}, \cite{Bi},
\cite{Comf}, \cite{Hew} and \cite{H}. Also, we must recall the
Bernstein's problem on approximation by polynomials on the whole
real line (see \cite{L}, \cite {L1} and \cite{L2}), and the
approximation problem for unbounded functions in $I$ (see for
example, \cite{ging}).

In recent years  it has arisen a new focus on the generalizations of
Weierstrass' approximation theorem,  which uses the weighted
approximation. More precisely, if $I$ is a compact interval, the
approximation problem is studied with the norm $L^\infty(I,\, w)$
defined by
\begin{equation}
\label{mm} \|f\|_{L^{\infty}(I,\,w)}:=\esup_{x\in I} |f(x)| w(x)\,,
\end{equation}
where $w$ is a weight, i.e., a non-negative measurable function and
the convention $0\cdot \infty=0$ is used. Observe that $(\ref{mm})$
is not the usual definition of the $L^\infty$ norm in the context of
measure theory, although it is the correct definition  when we work
with weights (see e.g. \cite{BO} and \cite{DMS}).

Considering weighted norms has been proved to be interesting mainly
because of two reasons: first, they allow to wider the set of
approximable functions (since the functions in $L^\infty(I,w)$ can
have singularities where the weight tends to zero); and, second, it
is possible to find functions which approximate a given function
$f$, whose qualitative behavior is similar to  qualitative behavior
of $f$ at those points where the weight tends to infinity.  The
reader  may find in  \cite{PQRT1}, \cite{PQRT2}, \cite{PQRT21} and
\cite{R1} recent and detailed  studies about such subject.

Another special kind of approximation problems arises when we
consider simultaneous approximation which includes derivatives of
certain functions; this is the case of versions of Weierstrass'
theorem in  weighted Sobolev spaces. With respect to recent
developments about this subject we refer to \cite{PQRT2} and
\cite{PQRT21}: In the first paper, under enough general conditions
concerning the vector weights (the so called type 1) defined in a
compact interval $I$, the authors characterize to the closure in the
weighted Sobolev space $W^{(k, \infty)} (I ,w)$ of the spaces of
polynomials, $k$-differentiable functions, and infinite
differentiable functions, respectively. In the second paper, the
reader may find shaper results for the case $k=1$.

In this paper we obtain a new result on Weierstrass' approximation
theorem with weights when considering approximation in Hilbert
spaces. We consider a real and separable Hilbert space $\GG$, a
compact interval $I\subset \RR$, the space of all the $\GG$-valued
essentially bounded functions $L^{\infty}_{\GG}(I)$, a weakly
measurable function $w:I\rightarrow \GG$, the space of all
$\GG$-valued continuous functions $C(I;\GG)$, the space of all the
$\GG$-valued functions $L^{\infty}_{\GG}(I,w)$, which are bounded
with respect to the norm defined by
\begin{equation}
\label{e2} \|f\|_{L^{\infty}_{\GG}(I,w)}:=\esup_{t\in I}
\|(fw)(t)\|_{\GG}\,.
\end{equation}

The paper is organized as follows. In Section 2 we provide some
notation, necessary preliminaries and auxiliary results which will
be often used throughout the text, we shall use standard notation or
it will be properly introduced whenever needed. In Section 3 we
present the main result about approximation in $L^{\infty}_{\GG}(I,
w)$.

\section{Preliminaries.}

In what follows, $I$ stands for any compact interval in $\RR$. By
$l^{2}(\RR)$ we denote the real linear space of all sequences
$\{x_{n}\}_{n\in\ZZ_{+}}$ with \hbox{
$\sum_{n=0}^{\infty}\left|x_{n}\right|^{2}<\infty$,} and $(\GG,
\langle\,\cdot\,,\,\cdot\, \rangle_{\GG})$ stands for a real and
separable Hilbert space with associated norm denoted by
$\|\,\cdot\,\|_{\GG}$.

It is well-known that every real and  separable Hilbert space $\GG$
is isomorphic either to $\RR^{n}$ for some $n\in\NN$ or to
$l^{2}(\RR)$. In each case, $\GG$ has structure of commutative
Banach algebra with the coordinatewise operations. In the first
case, we have commutative Banach algebra with identity and the
second case, this is commutative Banach algebra without identity.
The reader is referred to \cite{deu}, \cite{H1} or \cite{Y}  for
more details about these statements.

However, in many applications this isomorphism is not interesting:
for instance, we may be dealing with Hilbert spaces of, say,
analytic or differentiable functions, and in this  interplay between
the Hilbert space structure and the properties of individual
functions, the second can be  fruitful. Despite cases as the
previous,  to know such isomorphism is valuable, because it allows
to determine as far as the properties of Hilbert spaces can be
useful for us, and furthermore, we can just think in $l^{2}(\RR)$ or
$\RR^{n}$ when we want.

\subsection{On weighted spaces.}

A detailed discussion about properties of weighted spaces may be
found in \cite{Ch}, \cite{din}, \cite{KO}, \cite{NA1} or \cite{P2}.
We recall here some important tools and definitions which will be
used throughout this paper.

\begin{definition}
A  scalar weight $w$ is a measurable function $w:\RR\longrightarrow
[0,\infty]$. If $w$ is only defined in $A\subset \RR$, we set $w:=0$
in $\RR\setminus A$.
\end{definition}

\begin{definition}
Given a measurable set $A\subset \RR$ and a scalar weight $w$, we define
the space $L^\infty(A,w)$ as the space of equivalence classes of
measurable functions $f:A\longrightarrow \RR$ with respect to the
norm
$$
\|f\|_{L^\infty(A,w)}:=\esup_{x\in A} |f(x)|w(x)\,.
$$
\end{definition}

This space inherits some properties from the classical Lebesgue
space $L^{\infty}(A)$ and it allows us to study new functions, which
could not be in the classical $L^{\infty}(A)$ (see, for example
\cite{Ch}, \cite{Ru}). Another properties of $L^{\infty}(A,w)$ have
a strong relation with the nature of the weight $w$: in fact, if
$A=I$ and $w$ has multiplicative inverse, (i.e. there exists  a
weight $w^{-1}:I\longrightarrow \RR$, such that $w(t)w^{-1}(t)=1,
\quad \forall t\in I$) then, it is easy to see that
$L^{\infty}(I,w)$ and $L^{\infty}(I)$ are isomorphic, since the map
$\Psi_{w}:L^{\infty}(I,w)\rightarrow L^{\infty}(I)$ given by
$\Psi_{w} (f)=fw$  is a linear and bijective isometry, and
therefore, $\Psi_{w}$ is also homeomorphism, or equivalently, for
all $ Y\subseteq L^{\infty}(I,w)$, we have
$\Psi_{w}(\overline{Y})=\overline{\Psi_{w}(Y)}$,  where  we take
each closure with respect  to the norms $L^{\infty}(I,w)$ and
$L^{\infty}(I)$, respectively. Also, for all $ A\subseteq
L^{\infty}(I)$,
$\Psi^{-1}_{w}(\overline{A})=\overline{\Psi^{-1}_{w}(A)}$ and
$\Psi^{-1}_{w}=\Psi_{w^{-1}}$. Then  using  Weierstrass' theorem we
have,
\begin{eqnarray}
\label{ecc}
\Psi^{-1}_{w}(\overline{\pp})=\overline{\Psi^{-1}_{w}(\pp)}&=&\{f\in
L^{\infty}(I,w): fw\in C(I)\}.
\end{eqnarray}

Unfortunately,  the last equality in (\ref{ecc}) does not allow to
obtain information on  local behavior of the functions $f\in
L^{\infty}(I,w)$ which can be approximated. Furthermore,  if
\hbox{$f\in L^{\infty}(I,w)$,} then in general $fw$ is not
continuous function, since its continuity also depends of the
singularities of weight $w$ (see \cite{Ch}, \cite{L}, \cite{L1},
\cite{R1}, \cite{PQRT1}, \cite{PQRT21}).

The next definition presents the classification of the singularities
of a scalar  weight $w$  done in \cite{PQRT21} to show the results about
density of continuous functions in the space $L^{\infty}(\sop w , w
)$.

\begin{definition}
Given a scalar weight $w$ we say that $a\in \sop w$ is a {\it singularity}
of $w$ $($or {\it singular} for $w)$ if
$$
\eliminf_{x\in \sop w, \,x\to a} w(x)=0\,.
$$
We say that a singularity $a$ of $w$ is of {\it type} $1$ if
$\es_{x\to a} w(x)=0$.

We say that a singularity $a$ of $w$ is of {\it type} $2$ if \hbox{
$0<\elimsup_{x\to a} w(x)<\infty$.}

We say that a singularity $a$ of $w$ is of {\it type} $3$ if \hbox{
$\elimsup_{x\to a} w(x)=\infty$.}

We denote by $S$ and $S_i$ $(i=1,2,3)$ respectively, the set of
singularities of $w$ and the set of singularities of $w$ of type
$i$.

We say that $a\in S_i^+$ $($respectively $a\in S_i^-)$ if $a$
verifies the property in the definition of $S_i$ when we take the
limit as $x\to a^+$ $($respectively $x\to a^-)$. We define
$S^+:=S_1^+\cup S_2^+\cup S_3^+$ and $S^-:=S_1^-\cup S_2^-\cup
S_3^-$.
\end{definition}

\begin{definition}
Given a scalar weight $w$, we define the {\it right regular} and {\it left
regular} points of $w$, respectively, as
$$
R^+:= \big\{a\in\sop w:\, \eliminf_{x\in \sop w, \,x\to a^+} w(x)>0
\big\}\,,$$
$$
R^-:= \big\{a\in\sop w:\, \eliminf_{x\in \sop w, \,x\to a^-} w(x)>0
\big\}\,.
$$
\end{definition}

The following result was proved in  \cite{PQRT21} and it states a
characterization for the functions in $L^{\infty}(\sop w ,w)$ which
can be approximated by continuous functions in norm $L^{\infty}(\sop
w ,w)$ for every $w$.

\begin{teo}(Portilla et al. [\cite{PQRT21}, Theorem 1.2]).
\label{teo1} Let $w$ be any scalar weight and
$$H_{0} := \left\{
\begin{array}{l} f\in L^{\infty}(\sop w ,w): \,f\, \text{ is continuous to the right at every point of } R^+, \\
\qquad \qquad \qquad \quad f\, \text{ is continuous to the left at every point of } R^-, \\
\qquad \qquad \qquad \quad \text{for each } a\in S^+,\;\;
\es_{x\to a^+} |f(x)-f(a)|\,w(x)=0 \,,  \\
\qquad \qquad \qquad \quad  \text{for each } a\in S^-,\;\; \es_{x\to
a^-} |f(x)-f(a)|\,w(x)=0 \end{array}\right\}.
$$

Then:
\begin{enumerate}
\item[(a)] The closure of $C(\RR)\cap L^{\infty}(w)$ in
$L^{\infty}(w)$ is $H_0$.

\item[(b)] If $w\in L^\infty_{loc}(\RR)$, then the closure of
$C^\infty(\RR)\cap L^{\infty}(w)$ in $L^{\infty}(w)$ is also $H_0$.

\item[(c)] If $\sop w$ is compact and $w\in L^\infty(\RR)$, then
the closure of the space of polynomials is $H_0$ as well.

\end{enumerate}
\end{teo}

Theorem $\ref{teo1}$ is going to be  an important tool which is
going to allow us to obtain the key for the first result about
Hilbert extensions of Weierstrass' theorem with weights in the
present paper.

\subsection{ $\GG$-valued functions.}

\begin{definition}
Let $\GG$ be a real and  separable Hilbert space and we consider any
sequence  $\{x_{n}\}\subset \GG$. We say that the support of
$\{x_{n}\}$ is the set of $n$ for which $x_{n}\not=0$. We denote to
support of $\{x_{n}\}$ by $\sop{x_{n}}$.
\end{definition}

\begin{definition}
\label{def21} Let $\GG$ be a real and separable Hilbert space, a
weight $w$ on $\GG$  is  a weakly measurable function
$w:I\longrightarrow \GG$.
\end{definition}

Let $\GG$ be a real and separable Hilbert space. A $\GG$-valued
polynomial on $I$ is a function $\phi:I\rightarrow \GG$, such that
$$\phi(t)=\sum_{n\in\NN}\xi_{n}t^{n},$$
where $(\xi_{n})_{n\in\NN}\subset \GG$ has finite support.

Let $\pp(\GG)$ be the space of all $\GG$-valued polynomials on $I$.
It is well-known that $\pp(\GG)$ is a subalgebra of the space
$C(I;\GG)$ of all continuous $\GG$-valued functions on $I$.

For $1\leq p\leq\infty$, $L^{p}_{\GG}(I)$ denotes the set of all
weakly measurable functions $f:I\rightarrow \GG$ such that
$$
\int_{I}\|f(t)\|^{p}_{\GG}dt<\infty,\, \mbox{ if } \, 1\leq
p<\infty,
$$
or
$$\esup_{t\in I} \|f(t)\|_{\GG}<\infty,\, \mbox{ if } \, p=\infty.$$

Then $L^{2}_{\GG}(I)$ is a Hilbert space with respect to the inner
product
$$\langle f,g\rangle_{L^{2}_{\GG}(I)}=\int_{I}\langle f(t),g(t)\rangle_{\GG}\,dt.$$
$\pp(\GG)$ is also dense in $L^{p}_{\GG}(I)$, for $1\leq p< \infty$.

More details about these spaces may be found in \cite{na}.

\begin{definition}
\label{def21} Let $\GG$ be a real and separable Hilbert space, a
weight $w$ on $\GG$  is  a weakly measurable function
$w:I\longrightarrow \GG$.
\end{definition}

\begin{definition}
\label{def22} Let $w$ be a weight on $\GG$, we define the space
$L^{\infty}_{\GG}(I,w)$ as the space of equivalence classes of all
the  $\GG$-valued  weakly measurable functions $f:I\longrightarrow
\GG$ with respect to the norm
$$\|f\|_{L^{\infty}_{\GG}(I,w)}:=\esup_{t\in I} \|(fw)(t)\|_{\GG},\,$$
where $fw:I\longrightarrow \GG$ is defined as follows: If
$\dim\GG<\infty$; we have the functions $f$ and $w$ can be expressed
by $f=(f_{1},\ldots,f_{n_{0}})$ and $w=(w_{1},\ldots,w_{n_{0}})$,
respectively,  where $f_{j}, \,w_{j}:I\longrightarrow \RR$, for
$j=1,\ldots, n_{0}$, with $n_{0}=\dim \GG$. Then
$$(fw)(t):=
(f_{1}(t)w_{1}(t),\ldots,f_{n_{0}}(t)w_{n_{0}}(t)), \mbox{ for }
t\in I.$$

If $\dim\GG=\infty$, let $\{\tau_{j}\}_{j\in\ZZ_{+}}$ be a complete
orthonormal system, then for $t\in I$ the functions $f$ and $w$ can
be expressed as $f(t)=\sum_{j=0}^{\infty}\langle
f(t),\tau_{j}\rangle_{\GG}\, \tau_{j}$ and
$w(t)=\sum_{j=0}^{\infty}\langle w(t),\tau_{j}\rangle_{\GG}\,
\tau_{j},$ respectively. So, we can define
$$(fw)(t):=\sum_{j=0}^{\infty}\langle f(t),\tau_{j}\rangle_{\GG}\,\langle
w(t),\tau_{j}\rangle_{\GG}\, \tau_{j}, \mbox{ for } t\in I.$$
\end{definition}

On this way, we can  study our approximation problem using the
properties of commutative Banach algebra of $l^{2}(\RR)$.

The next Proposition shows a result about algebraic properties and
density of $\pp(\GG)$ in $C(I;\GG)$. The analogous result,
when $\GG$ is a  complex separable  Hilbert space,  appears in
\cite{na}.

\begin{prop}
\label{pro1} \hspace{1cm}
\begin{enumerate}
\item[i)]$\pp(\GG)$ is a subalgebra of the space of all
$\GG$-valued continuous functions on $I$. \item[ii)] The closure of
$\pp(\GG)$ in $L^{\infty}_{\GG}(I)$ is $C(I;\GG)$.
\end{enumerate}
\end{prop}
\begin{proof}
\hspace{1cm}
\begin{enumerate}
\item[i)] It is straight forward.

\item[ii)] It is enough to prove that
$C(I;\GG)\subset\overline{\pp(\GG)}$, since
$\overline{\pp(\GG)}\subset \overline{C(I;\GG)}=C(I;\GG)$.

{\bf Case 1:} $\dim\GG<\infty$.

Let us assume that $\dim\GG=n_{0}$. Given  an orthonormal basis $\{\tau_{1},\ldots,
\tau_{n_{0}}\}$ of $\GG$, $\epsilon>0$ and
$f\in C(I;\GG)=C(I,\RR^{n_{0}})$, then
$f\sim(f_{1},\ldots,f_{n_{0}})$ with $f_{j}\in C(I)$,
$j=1,\ldots,n_0$. The Weierstrass' theorem guarantees that there
exists  $p_{k}\in \pp$ such that
$$\|f_{j}-p_{j}\|_{L^{\infty}(I)}<\frac{\epsilon}{\sqrt{n_{0}}}, \quad j=1,\ldots, n_{0}.$$

If we consider the polynomial $p\in\pp(\GG)$ such that  $p\sim
(p_{1},\ldots,p_{n_{0}})$, then we have that
\begin{eqnarray*}
\|f-p\|_{L^{\infty}_{\GG}(I)}&=&\esup_{t\in I}\|(f-p)(t)\|_{\GG}\\
&=& \esup_{t\in I} \left[\sum_{j=1}^{n_{0}}|\langle
f(t)-p(t),\tau_{j} \rangle_{\GG}|^{2}
\right]^{1/2}\\
&\leq& \esup_{t\in I}\|((f_{1}-p_{1})(t),\ldots,
(f_{n_{0}}-p_{n_{0}})(t))\|_{\RR^{n_{0}}}<\epsilon.
\end{eqnarray*}

{\bf Case 2:} $\GG$ is infinite-dimensional.

Let $f\in C(I;\GG)$ and $\{\tau_{j}\}_{j\in\ZZ_{+}}$ a  complete
orthonormal system, then for each $t\in I$
$$f(t)=\sum_{j=0}^{\infty}\langle f(t),\tau_{j}\rangle_{\GG}\, \tau_{j},$$
consequently, given $\epsilon>0$ there exists $m_0 \in \ZZ_{+}$ such
that
$$\left\|f(t)-\sum_{j=0}^{n}\langle f(t),\tau_{j}\rangle\tau_{j}\right\|_{\GG}<\epsilon, \mbox{ whenever } n\geq m_{0}.$$

Now, let us consider the functions $f_{j}:I\rightarrow\RR$ given by
$f_{j}(t)=\langle f(t), \tau_{j}\rangle_{\GG}$. We have that
$f\sim\{f_{j}\}$ with  $\sum_{j\in\ZZ_{+}}|f_{j}(t)|^{2}<\infty$,
for each $t\in I$ and $f_{j}\in C(I)$.

So, Weierstrass' approximation theorem guarantees that there exists
a sequence $\{p_{j}\}_{j\in\ZZ_{+}}\subset \pp$ such that
$$\|f_{j}-p_{j}\|_{L^{\infty}(I)}<\frac{\epsilon}{j+1}, \quad j\in\ZZ_{+}.$$

We define the $\GG$-polynomials $\tilde{p}_{j}\in \pp(\GG)$ by
$\tilde{p}_{j}(t)=p_{j}(t)\tau_{j}$, for each $j\in\ZZ_{+}$. Then
for $n\geq m_{0}$ we have

\begin{eqnarray*}
\left\|f(t)-\sum_{j=0}^{n}\tilde{p}_{j}(t)\right\|_{\GG}&\leq& \left\|f(t)-\sum_{j=0}^{n}f_{j}(t)\tau_{j}\right\|_{\GG}+\left\|\sum_{j=0}^{n}f_{j}(t)\tau_{j}-\sum_{j=0}^{n}\tilde{p}_{j}(t)\right\|_{\GG}\\
&\leq& \epsilon +\left( \sum_{j=0}^{\infty}\left|f_{j}(t)-p_{j}(t)\right|^{2}\right)^{1/2}\\
&<&\epsilon
+\left(\sum_{j=0}^{\infty}\left(\frac{\epsilon}{j+1}\right)^{2}\right)^{1/2}=\epsilon\left(
1+\left(\sum_{j=0}^{\infty}\frac{1}{(j+1)^{2}}\right)^{1/2}\right).
\end{eqnarray*}

From these inequalities we can deduce that for $n$ large enough
there exists $q_{n}(t)=\sum_{j=0}^{n}\tilde{p}_{j}(t)$ such that
$$\|f-q_{n}\|_{L^{\infty}_{\GG}(I)}<C\epsilon.$$
\end{enumerate}

This completes the proof.
\end{proof}

\section{The main results.}

In this section, we only deal with weights $w$ such that $\sop w=I$.

\subsection{Approximation in $L_{\GG}^{\infty}(I,w)$.}
\begin{definition}
Let $\GG$ be  a real and  separable Hilbert space and let $w$ be a
weight on $\GG$. We say that $w$ is admissible$^{\star}$ if one of
the following conditions is satisfied
\begin{enumerate}
\item[i)] If $\dim\GG<\infty$ then each one of the components
$w_j$, $1\leq j\leq \dim\GG$, is a scalar weight.
\item[ii)] If $\dim\GG=\infty$, $\{\tau_{j}\}_{j\in\ZZ_{+}}$ is a complete
orthonormal system, and $w(t)=\sum_{j=0}^{\infty}\langle
w(t),\tau_{j}\rangle_{\GG}\, \tau_{j}$, then each one of the
functions $\langle w(t),\tau_{j}\rangle_{\GG}$ is a weight function.
\end{enumerate}
\end{definition}

Let us observe that if $\dim\GG=\infty$ and $w$ is
admissible$^{\star}$, then it induces a family of weighted
$l^{2}(\RR)$ spaces, $\{l^{2}_{t}(\RR;w): t\in I\}$ given by
$$l^{2}_{t}(\RR;w)=\left\{\{x_{j}\}_{j\in\ZZ_{+}}: \sum_{j=0}^{\infty} \langle
w(t),\tau_{j}\rangle_{\GG}\,|x_{j}|^{2}<\infty \right\}.$$

And for each $t\in I$, the function $w_{j}(t)=\langle
w(t),\tau_{j}\rangle_{\GG}$ also induces a linear isometry

$\Psi^{t}_{w_{j}}:l^{2}_{t}(\RR;w_{j})\rightarrow l^{2}(\RR)$ given
by
$$\Psi^{t}_{w_{j}}\left( \{x_{j}\}_{j\in\ZZ_{+}}\right)= \{w_{j}(t)x_{j}\}_{j\in\ZZ_{+}}=
\{\langle w(t),\tau_{j}\rangle_{\GG}\,x_{j}\}_{j\in\ZZ_{+}}.$$

For a brief study on weighted $l^{2}(\RR)$ spaces the reader is
referred to \cite{deu}. In order to characterize the $\GG$-valued
functions which can be approximated in $L_{\GG}^{\infty}(I,w)$ by
functions in $C(I;\GG)\cap L_{\GG}^{\infty}(I,w)$, our argument
requires an admissible$^{\star}$ weight  $w$. It is clear that in
the one-dimensional case an admissible$^{\star}$ weight is an
arbitrary scalar weight on $I$, and therefore the Theorem 2.1 in
\cite{PQRT21} holds in this case.

\begin{teo}
\label{main1} Let $\GG$ be a real and  separable Hilbert space and
let $w$ be an admissible$^{\star}$ weight on $\GG$. Let us define
$$H:=\left\{\begin{array}{l} f\in L^{\infty}_{\GG}(I,w):
f\sim (f_{1},\ldots,f_{n_{0}})\, \mbox{ and }\, f_{j}\in H_{j},\,
1\leq j\leq n_{0}\, \mbox{ with }\, n_{0}=\dim\GG,\\
\qquad \qquad \qquad \mbox{ or }\, f\sim \{f_{j}\}\, \mbox{ and }\,
f_{j}\in H_{j},\, j\in\ZZ_{+}\,  \mbox{ if }\,
\dim\GG=\infty.\end{array}\right\},$$ where
$$H_{j} := \left\{
\begin{array}{l} f_{j}\in L^{\infty}(I,w_{j}): \,f_{j}\, \text{ is continuous to the right at every point of } R^+, \\
 \qquad \qquad \qquad \quad f_{j}\, \text{ is continuous to the left at every point of } R^-, \\
 \qquad \qquad \qquad \quad \text{for each } a\in S^+,\;\;
\es_{x\to a^+} |f_{j}(x)-f_{j}(a)|\,w_{j}(x)=0 \,,  \\
\qquad \qquad \qquad \quad  \text{for each } a\in S^-,\;\; \es_{x\to
a^-} |f_{j}(x)-f_{j}(a)|\,w_{j}(x)=0 \end{array}\right\}.
$$
Then the closure of $C(I;\GG)\cap L^{\infty}_{\GG}(I,w)$ in
$L^{\infty}_{\GG}(I,w)$ is $H$. Furthermore, if $w\in
L_{\GG}^{\infty}(I)$ then the closure of the space of $\GG$-valued
polynomials is $H$ as well.
\end{teo}

\begin{proof} Let us assume that $\dim\GG=n_{0}$. If $f\in \overline{C(I;\GG)\cap
L^{\infty}_{\GG}(I,w)}^{L^{\infty}_{\GG}(I,w)}$, then $f\sim
(f_{1},\ldots, f_{n_{0}})$, with $f_{j}:I\longrightarrow\RR$, $1\leq
j\leq n_0$. Given $\epsilon
>0$, there exists $g\in C(I;\GG)\cap
L^{\infty}_{\GG}(I,w)$ such that
$\|f-g\|_{L^{\infty}_{\GG}(I,w)}<\epsilon$. Let us consider
$(g_{1},\ldots,g_{n_{0}})$ such that $g_{j}\in C(I)\cap
L^{\infty}(I,w_{j})$ and $g\sim (g_{1},\ldots,g_{n_{0}})$, then
$$|(f_{j}(t)-g_{j}(t))w_{j}(t)|\leq \esup_{s\in
I}\left[\sum_{j=1}^{n_{0}}|(f_{j}(s)-g_{j}(s))w_{j}(s)|^{2}
\right]^{1/2} \, \mbox{ a.e. }$$

On other hand,
$\left[\sum_{j=1}^{n_{0}}|(f_{j}(s)-g_{j}(s))w_{j}(s)|^{2}
\right]^{1/2}=\|f-g\|_{L^{\infty}_{\GG}(I,w)}$, as consequence of
$\GG$ is isomorphic to $\RR^{n_{0}}$ and the Parseval identity (see
\cite{deu} or \cite{Y}). Then,
$$\|f_{j}-g_{j}\|_{L^{\infty}(I,w_{j})}\leq
\|f-g\|_{L^{\infty}_{\GG}(I,w)}<\epsilon.$$

Hence, $f_j\in \overline{C(I)\cap L^{\infty}(I,w_{j})}^{
L^{\infty}(I,w_{j})}$ for $1\leq j\leq n_{0}$, and the \hbox{part
(a)} of Theorem \ref{teo1} gives that $H$ contains
$\overline{C(I;\GG)\cap
L^{\infty}_{\GG}(I,w)}^{L^{\infty}_{\GG}(I,w)}$.

In order to see that $H$ is contained in $\overline{C(I;\GG)\cap
L^{\infty}_{\GG}(I,w)}^{L^{\infty}_{\GG}(I,w)}$, let us fix $f\in H$
and $\epsilon>0$, and let us consider each one of its component functions $f_j \in
H_j$, $j=1,\ldots, n_{0}$. By the part (a) of Theorem \ref{teo1},
there exists $g_j \in C(I)\cap L^{\infty}(I,w_{j})$, $j=1,\ldots,
n_{0}$,  such that
$$\|f_{j}-g_{j}\|_{L^{\infty}(I,w_{j})}<\frac{\epsilon}{\sqrt{n_{0}}}.$$
We define the function $g\in C(I;\GG)$  such that  $g\sim
(g_{1},\ldots,g_{n_{0}})$, then
\begin{eqnarray*}
\|f-g\|_{L^{\infty}_{\GG}(I,w)}&=&\esup_{t\in I}\|((f-p)w)(t)\|_{\GG}\\
&=&  \esup_{t\in
I}\left[\sum_{j=1}^{n_{0}}|(f_{j}(t)-g_{j}(t))w_{j}(t)|^{2}
\right]^{1/2}<\epsilon.
\end{eqnarray*}

If $w\in L^{\infty}_{\GG}(I)$, the closure of the $\GG$-valued
polynomials is $H$ as well, as a consequence of  Proposition
\ref{pro1}.

In a similar way, if $\dim\GG=\infty$, $\{\tau_{j}\}_{j\in\ZZ_{+}}$ is
a complete orthonormal system and $f\in \overline{C(I;\GG)\cap
L^{\infty}_{\GG}(I,w)}^{L^{\infty}_{\GG}(I,w)}$, then $f(t)=
\sum_{j=0}^{\infty}\langle f(t),\tau_{j}\rangle \tau_{j}$. Given
$\epsilon
>0$, there exists $g\in C(I;\GG)\cap
L^{\infty}_{\GG}(I,w)$ such that \hbox{
$\|f-g\|_{L^{\infty}_{\GG}(I,w)}<\epsilon$.} Let us consider
$\{g_{j}\}_{j\in\ZZ_{+}}$ such that $g_{j}\in C(I)\cap
L^{\infty}(I,w_{j})$ and \hbox{$g\sim \{g_{j}\}_{j\in\ZZ_{+}}$,}
then
$$|(f_{j}(t)-g_{j}(t))w_{j}(t)|\leq \esup_{s\in
I}\left[\sum_{j=0}^{\infty}|(f_{j}(s)-g_{j}(s))w_{j}(s)|^{2}
\right]^{1/2} \, \mbox{ a.e. }$$

On other hand,
$\left[\sum_{j=1}^{\infty}|(f_{j}(s)-g_{j}(s))w_{j}(s)|^{2}
\right]^{1/2}=\|f-g\|_{L^{\infty}_{\GG}(I,w)}$, as consequence of
$\GG$ is isomorphic to $l^{2}(\RR)$ and the Parseval identity (see
\cite{deu} or \cite{Y}). Then,
$$\|(f_{j}-g_{j}\|_{L^{\infty}(I,w_{j})}\leq
\|f-g\|_{L^{\infty}_{\GG}(I,w)}<\epsilon.$$

Hence, $f_j\in \overline{C(I)\cap L^{\infty}(I,w_{j})}^{
L^{\infty}(I,w_{j})}$ for $j\in\ZZ_{+}$, and  the part (a) of
Theorem \ref{teo1} gives that $H$ contains $\overline{C(I;\GG)\cap
L^{\infty}_{\GG}(I,w)}^{L^{\infty}_{\GG}(I,w)}$.

In order to see that $H$ is contained in $\overline{C(I;\GG)\cap
L^{\infty}_{\GG}(I,w)}^{L^{\infty}_{\GG}(I,w)}$, let $f\in H$ and
$\epsilon>0$, and let us consider the component functions $f_j \in
H_j$ of $f$, $0\leq j<\infty$. Since $w_{j}(t)=\langle
w(t),\tau_{j}\rangle$ is a weight, by the part (a) of Theorem
\ref{teo1}, there exists $g_j \in C(I)\cap L^{\infty}(I,w_{j})$,
$0\leq j<\infty$,  such that
$$\|f_{j}-g_{j}\|_{L^{\infty}(I,w_{j})}<\frac{\epsilon}{j+1}, \,\,\quad  j\in \ZZ_{+}.$$
We define the function $g:I\longrightarrow \GG$ by
$g(t)=\sum_{j=0}^{\infty}g_{j}(t)\tau_{j}$, then
\begin{eqnarray*}
\|f-g\|_{L^{\infty}_{\GG}(I,w)}&=&\esup_{t\in I}
\|((f-g)w)(t)\|_{\GG}\\
&=&\esup_{t\in
I}\|\{(f_{j}(t)-g_{j}(t))w_{j}(t)\}\|_ {l^{2}(\RR)}\\
&=&\esup_{t\in I}
\left[\sum_{j=0}^{\infty}|f_{j}(t)-g_{j}(t)|^{2}w_{j}^{2}(t)\right]^{1/2}\\
  &\leq& \left[\sum_{j=0}^{\infty}\left(\frac{\epsilon}{j+1}\right)^{2}\right]^{1/2}
  =\epsilon\left[\sum_{j=0}^{\infty}
\frac{1}{(j+1)^{2}}\right]^{1/2}
\end{eqnarray*}
\end{proof}

This result is similar when $\GG$ is a complex separable Hilbert
space and it can also be extended to $L^{\infty}_{L(\GG)}(I,w)$,
where $L(\GG)$ is the space of operators on $\GG$.

{\bf Acknowledgment.}\\
The author wishes to thank  the referee and Professor Jos\'e Manuel
Rodr\'{\i}guez for their suggestions and comments which have
improved the presentation of the paper.

\vspace{1cm}
\begin{tabular}{c}
\small
\parbox{7cm}{
\quad\\
\quad\\
Yamilet Quintana\\
Departamento de Matem\'aticas\\
Puras y Aplicadas\\
Edificio Matem\'aticas y Sistemas (MYS)\\
Apartado Postal: 89000, Caracas 1080 A\\
Universidad Sim\'on Bol\'{\i}var\\
Caracas\\
VENEZUELA\\
\quad\\}
\end{tabular}
\end{document}